\documentclass[11pt,twoside,reqno]{article}
\usepackage{latexsym,epsfig,graphpap,graphics,amsmath,amssymb,float}
\usepackage{hyperref}
\usepackage[utf8x]{inputenc}
\newtheorem{thm}{Theorem}[section]

\newtheorem{lem}[thm]{Lemma}
\newtheorem{prop}[thm]{Proposition}
\newtheorem{defn}[thm]{Definition}
\newtheorem{remark}[thm]{Remark}

\newtheorem{example}[thm]{Example}
\numberwithin{equation}{section}

\usepackage{plain}
\begin{document}
\title{Approximations for Pareto and Proper Pareto solutions and their KKT conditions}
\author{P. Kesarwani\footnote{Department of Mathematics and Statistics, Indian Institute Of Technology Kanpur,India}, P.K.Shukla \footnote{Institute AIFB, Karlsruhe Institute of Technology, Karlsruhe, Germany}, J. Dutta\footnote{Department of Economics Sciences, Indian Institute Of Technology Kanpur, India} and K. Deb\footnote{College of Engineering, Michigan State University, Michigan, USA}}

\maketitle

\begin{abstract}
In this article, we view the approximate version of Pareto and weak Pareto solutions of the multiobjective optimization problem through the lens of KKT type conditions. We also focus on an improved version of Geoffrion proper Pareto solutions and characterize them through saddle point and KKT type conditions. We present an approximate version of the improved Geoffrion proper solutions and propose our results in general settings.
\end{abstract}

\section{Introduction}
The importance of multiobjective optimization problems in various applications in engineering, business and management can be hardly overstated. For a wide range of applications in engineering design see, for example, \cite{deb2001multi}.

\noindent From a theoretical point of view, the idea of multiobjective optimization becomes challenging since we are speaking about minimizing/maximizing a vector-valued function. In order to define the notion of a solution, we need to depend on the partial order, which is often induced by a closed convex pointed cone on the image space of the objective function. This leads to two fundamental notions of solutions, namely the Pareto solutions and weak Pareto solutions. The points corresponding to these solutions in the image space of the objective function are often referred to as efficient solutions and weak efficient solutions. The collection of all efficient solutions is often referred to as the Pareto efficient frontier. We emphasize that the notion of Pareto and weak Pareto solutions are global notions. Mathematically speaking, it is not at all difficult to devise a local counterpart, and it is the global aspect that is sought by the decision makers. Further, the idea of Pareto solutions is often considered more relevant than the weak solutions from the point of view of the applications. We refer the following monographs of Ehrgott \cite{EHR}, Jahn \cite{JAHN}, Luc \cite{luc1987scalarization}, Chankong et al. \cite{chankong2008multiobjective} and the references therein to see the development of multiobjective optimization over the past several decades.\\
There are several approaches to solve multiobjective problems, for example, scalarization methods, descent methods, metaheuristics and many more. But when it comes to actual computation using the mentioned methods, the algorithms always produce approximate solutions. Thus, it is essential to define notions of approximate solutions and characterize their properties. There are various notion of approximate solutions in the literature (see \cite{loridan1984varepsilon}, \cite{valyi1985approximate}, \cite{dutta2001approximate}, \cite{gutierrez2006approximate}, \cite{gutierrez2010optimality}) which deals with characterizing introduced notions in greater details. In this article, our main aim is to revisit the fundamental notion of approximate Pareto solutions and a proper Pareto solution (which we shall describe below). Further, we analyze these solutions through KKT type conditions. This approach to studying KKT type conditions for approximate solutions can lead to the development of stopping criteria for algorithms. It can also be used to check the quality of the approximate solution produced by any algorithm used to solve the multiobjective problem. \\
It is essential to a decision maker, who is taking some decisions based on multiobjective optimization models need not necessarily be interested in all the Pareto solutions of the problem at hand.
In many cases, the decision maker focuses on the part of the Pareto frontier in the image space, which corresponds to a subset of the set of Pareto solutions. These subsets, when chosen in a particular way, gives rise to various classes of proper Pareto solutions (see\cite{EHR}). Very recently, the authors discussed an improved version of Geoffrion proper solutions in \cite{shukla2019practical}. This solution notion is based on the assumption that the decision maker, in practice, usually looks for those proper solutions whose trade-off bound is bounded by a value preset by her/him. The detailed analysis of such solutions and their approximate version has been carried out and shown to be stable than the standard Geoffrion solutions (for more details see \cite{shukla2019practical}). In the present article, our major goal is to analyse saddle point and KKT type conditions for these solutions.\\
The whole paper revolves around answering three questions in which first two questions stem from an attempt to generalize two results, which are on approximate solutions for scalar optimization problems which appeared in \cite{dutta2013approximate}. The first result concerns a scalar optimization problem with locally Lipschitz data (see Theorem 3.2 in \cite{dutta2013approximate}) which says that if a sequence of points each satisfying an approximate version of KKT conditions converges to point under a suitable constraint qualification, then the limit of the sequence is a KKT point. Thus, we have the following first question:
\begin{itemize}

\item {\textbf{Q1:} Can a similar kind of result be deduced for multiobjective optimization problem?}

\end{itemize}
Our second question stems from Theorem 3.7 in \cite{dutta2013approximate} in which the reversed result of the Theorem 3.2 is asked. The result conclude an affirmative answer for the reverse result which proves that for any local minimizer  of an optimization problem with suitable constraint qualification, there a sequence of points converge to that local minima and there exists a subsequence of the main sequence which satisfies some type of approximate KKT type conditions when they are very near to the solution.

\begin{itemize}

\item {\textbf{Q2:} Can we generalize the reverse result in the multiobjective settings? Further, do the locally Lipschitz data suffice, or we need more assumptions? Can the convexity assumption give us better results? }

\end{itemize}
Our third question is associated with the KKT-type conditions for the approximate Geoffrion proper solutions with a preset bound.

\begin{itemize}

\item {\textbf{Q3:} Can we develope an approximate KKT type condition which can completely characterize a Geoffrion proper solutions with a preset bound at least in the convex case? Does the saddle point conditions completely characterize such class of solutions? }
\end{itemize}

The paper is organised as follows. In Section 2, we present the problem, basic definitions and the technical tools from convex and non-smooth analysis required in the article. In Section 3, we answer the first two questions raised in this section and Section 4 assures the last question by trying to develop the saddle point conditions and approximate KKT type conditions for the improved Geoffrion proper solutions. We end our discussion by concluding remarks in Section 5. We want to end this section by stating that most symbols used in the article are fairly standard in the literature.
\section{Preliminaries and basic tools}
Let $A\subseteq \mathbb{R}^n$ be a given set, then closure and interior of
set $A$ is denoted by cl$A$ and int$A$ respectively. For vectors $x,y \in \mathbb{R}^n$ the inner product given by $\langle x, y\rangle$. A set $A\subset \mathbb{R}^n$ is a cone, if for each $a\in A$ and positive scalar $\lambda$, $\lambda a \in A$. A cone $A$ is pointed, if $A\cap (-A)=\{0\}$. A normal cone of a convex set $A$ at the point $x_0$, denoted by $N_A(x_0)$, is $N_A(x_0)=\{v\in \mathbb{R}^n: \langle v, x-x_0\rangle \leq 0, \text{ for all }  x\in U \}$.
We consider the following form of multiobjective optimization problem (MOP) in this article:
\begin{eqnarray*}
&&\min f(x):=(f_1(x),\ldots,f_m(x)),\\
&&{\rm subject~to}~~ g_j(x)\leq 0,~ j=1,2,..,l.
\end{eqnarray*}
where each $f_i:\mathbb{R}^n\to \mathbb{R}$ and $g_j:\mathbb{R}^n\to \mathbb{R}$. Let us denote the constraint set by $X:=\{x\in \mathbb{R}^n: g_j(x)\leq 0,~j=1,2,..,l\}\subseteq \mathbb{R}^n$, $I:=\{1,2,..,m\}$, $L:=\{1,2,..,l\}$. As we mentioned earlier that there are several notions for approximate solutions but in this article, we consider the notion of approximate solution introduced in Loridan \cite{loridan1984varepsilon}. We consider $\epsilon\in\mathbb{R}^m_+,$ \textit{i.e.}, $\epsilon=(\epsilon_1,\epsilon_2,\ldots,\epsilon_m)$, $\epsilon_i\geq 0$ for each $i\in I$ to formalize our notions. Our focus on this paper is on $\epsilon$- solutions of MOP. The partial order of image space $f(X)\subseteq \mathbb{R}^m$ is induced by natural cone $\mathbb{R}^m_+$ in the following definition.  

\begin{defn}\label{d3}
Given $\epsilon\in \mathbb{R}^m_+$, if there is no $x\in X$ such that $f(x)+\epsilon-f(x^*)\in -\mathbb{R}^m_+\setminus\{0\},$ then the point $x^*\in X$ is said to be an $\epsilon$-Pareto optimal solution of MOP. Further if there is no $x\in X$ such that
$f(x)+\epsilon-f(x^*)\in -{\rm int}(\mathbb{R}^m_+),$
then the point $x^*$ is said to be a weak $\epsilon$-Pareto optimal solution of MOP. 
\end{defn}
\noindent An $\epsilon$-Pareto (weak) optimal solution with $\epsilon=0$ is commonly known as Pareto (weak) optimal solution. Though not always seen in the literature the following notions of a local solutions are also relevant.
\begin{defn}\label{d2}
A point $x^*$ is said to be a loacl Pareto optimal solution of MOP if there exists $\delta >0$ and no $x\in X\cap B_{\delta}(x^*) $ such that,
$f(x)-f(x^*)\in -\mathbb{R}^m_+\setminus\{0\},$
where $B_{\delta}(x_0)\subset \mathbb{R}^n$ is a ball of radius $\delta$.
\end{defn}
\noindent The weak counter part of local solution can be defined in the similar fashion as in Definition \ref{d2}.
We want to mention that in several situations we consider the particular form of the vector $\epsilon \in \mathbb{R}^m_+$, given by $\epsilon=\varepsilon e$, where $e=(1,1..,1)^T$ and $\varepsilon\in \mathbb{R}_+$. In those cases, the solutions referred to as the $\varepsilon e$-Pareto and $\varepsilon e$-weak Pareto solution respectively. The set of all $\epsilon$-Pareto points is denoted by ${S}_{\epsilon}(f, X)$ and the set of all $\epsilon$-weak Pareto points as $S_{w,\epsilon}(f, X)$.
 

\begin{defn}\label{d7}
Given $\epsilon\in \mathbb{R}^n_+$, a point $x_0\in X$ is called $\epsilon$-Geoffrion proper solution of MOP if $x_0\in \mathcal{S}_{\epsilon}(f,X)$  and if there exists a number $M>0$ such that for all $i\in I$ and $x\in X$ satisfying $f_i(x)< f_i(x_0)-\epsilon_i$, there exists an index $j\in I$ such that $f_j(x_0)-\epsilon_j< f_j(x)$ and
\begin{eqnarray*}
\frac{f_i(x_0)- f_i(x)-\epsilon_i}{f_j(x)- f_j(x_0)+\epsilon_j}\leq M.
\end{eqnarray*}\end{defn}

\noindent The upper bound of the trade-off in the above definition is not known beforehand and the definition only assures the existence of such a bound. Further, it is clear form the definition that the trade-off varies as we choose different proper points. The improved definition introduced in \cite{shukla2019practical} eliminates the dependence of the bound on the solution points. Let us state the improved notion of Geoffrion proper solutions studied in \cite{shukla2019practical}.

\begin{defn}\label{d13}
Given $\epsilon\in \mathbb{R}^n_+$ and a scalar $\hat{M}>0$, a point $x_0\in X$ is called $(\hat{M},\epsilon)$-Geoffrion  proper solution of MOP if $x_0\in \mathcal{S}_{\epsilon}(f,X)$  and for all $i\in I$ and $x\in X$ satisfying $f_i(x)< f_i(x_0)-\epsilon_i$, there exists an index $j\in I$ such that $f_j(x_0)-\epsilon_j< f_j(x)$ and
\begin{eqnarray*}
\frac{f_i(x_0)- f_i(x)-\epsilon_i}{f_j(x)- f_j(x_0)+\epsilon_j}\leq \hat{M}.
\end{eqnarray*}\end{defn}
Given $\hat{M}>0$, we shall denote the set of all $(\hat{M},\epsilon)$- Geoffrion proper as $\mathcal{G}_{\hat{M},\epsilon}(f,X)$. For $\epsilon=0$, the set of exact $\hat{M}$- Geoffrion proper is denoted by $\mathcal{G}_{\hat{M}}(f,X)$. Now we shall present the Ekeland variation principle for vector-valued functions which was introduced in \cite{tammer1992generalization} when the ordering cone is $\mathbb{R}_+^m$. We first define the notion of lower semicontinuity and boundness of vector-valued functions which will be needed in the principle.

\begin{defn}\label{d1}
Let $f:U\rightarrow \mathbb{R}^m$ where $U$ is a non-empty subset of $\mathbb{R}^n$. The function $f$ is {$\mathbb{R}^m_+$-bounded below} if there exists $y\in \mathbb{R}^m$ such that $f(x)-y\in \mathbb{R}_+^m$ for all $x\in U$. Let $c\in int(\mathbb{R}_+^m)$, the function $f$ is {$(c,\mathbb{R}_+^m)$- lower semi continuous} if for all $t\in \mathbb{R}$, $\{x\in U:tc- f(x)\in \mathbb{R}_+^m\}$ is closed.
\end{defn}
\begin{thm}\label{d17}
Let $f:U\rightarrow \mathbb{R}^m$ where $U\subseteq \mathbb{R}^m$ be a $(c_0,\mathbb{R}^m_+)$- lower semi continuous function for $c_0\in int(\mathbb{R}_+^m)$ which is also $\mathbb{R}^m_+$-bounded below . Further, suppose we are given 
  $\rho >0$ and a point $x_0\in U$ such that,
\begin{equation}\label{eqm}
f(x)+\rho c_0-f(x_0)\not \in -\mathbb{R}^m_+\setminus \{0\}, \text{ for all } x \in U.
\end{equation}
Then, there exists $\bar{x_0}=\bar{x_0}(\rho)\in U$ such that  $\|\bar{x_0}-x_0\|\leq \sqrt{\rho}$ and for all $x\in U\setminus \{\bar{x_0}\}$
\begin{enumerate}
\item  $f(x)+ \rho c_0-f(\bar{x_0})\not \in -int (\mathbb{R}^m_+) $,
\item  $f(x)+\sqrt{\rho}\|\bar{x_0}-x\| c_0-f(\bar{x_0})\not \in -int (\mathbb{R}^m_+)$. 
\end{enumerate}
\end{thm}
In this article, we rely on two major tools from non-smooth analysis, namely the subdifferential of a convex function and the Clarke subdifferential of a locally Lipschitz function. Though these notions are very well known in the optimization community, we shall provide the definitions for completeness. We shall however restrict ourselves to the class of functions which are finite-valued function on $\mathbb{R}^n$. \\
Let $f:\mathbb{R}^n\rightarrow \mathbb{R}$ be a convex function, then the subdifferential of $f$ at the point $x$ is a set of vectors in $\mathbb{R}^n$, given as 
$$\partial f(x)=\{v\in \mathbb{R}^n: f(y)-f(x)\geq \langle v,y-x\rangle,~\text{for all}~y\in \mathbb{R}^n\}.$$
The subdifferential set is a non-empty, convex and compact for every $x\in \mathbb{R}^n$. The subdifferential is also deeply linked with the notion of the directional derivative of a convex function. The directional derivative of a convex function at a given $x$ in the direction $h$ is given as 
$$f'(x,h)=\lim_{\lambda\downarrow 0} \frac{f(x+\lambda h)-f(x)}{\lambda}$$
This directional derivative exists for each $x$ and in each direction $h$, and, the subdifferential of $f$ can be written as $\partial f(x)=\{v\in \mathbb{R}^n:~f'(x,h)\geq \langle v,h \rangle,~\text{for all}~h\in \mathbb{R}^n\}.$ Thus each of these can be recovered from the other. The generalized notion of derivative has properties like the usual derivative of calculus. We will begin with the most fundamental one, the sum rule. Let $f:\mathbb{R}^n\rightarrow \mathbb{R}$ and $g:\mathbb{R}^n\rightarrow \mathbb{R}$ are convex functions. Then 
\begin{equation}\label{eq21}
\partial (f+g)(x)=\partial f(x)+\partial g(x).
\end{equation}
\noindent For more details on subdifferentail of convex functions see \cite{bazaraa2013nonlinear}. It is important to note that a point $x_0$ is a global minimum of $f$ on $\mathbb{R}^n$ if and only if $0\in \partial f(x_0)$. Since subdifferential is a generalized version of derivative, it has some limitation. The $\varepsilon$-subdifferential is a relaxed version of the subdifferential which is very useful tool in convex analysis and optimization. We begin with defining the $\varepsilon$-subdifferential of convex function.
\begin{defn}\label{d8}
Let $f:\mathbb{R}^n\rightarrow \mathbb{R}$ be convex function and $\varepsilon \geq 0$. The $\varepsilon$-subdifferential of $f$ at the point $x$ is given as 
$$\partial_{\varepsilon} f(x)=\{v\in \mathbb{R}^n: f(y)-f(x)\geq \langle v,y-x\rangle-\varepsilon,~\text{for all}~y\in \mathbb{R}^n\}.$$
\end{defn} 
The elements of $\partial_{\varepsilon} f(x)$ are called $\varepsilon$-gradients of $f$ at $x$ and $\partial_{\varepsilon} f(x)\not =\emptyset$ for all $x\in \mathbb{R}^n$. A point $x_0$ is called an $\varepsilon$-minimizer of $f$ on $\mathbb{R}^n$ if $f(y)-f(x)\geq -\varepsilon$, for all $y\in \mathbb{R}^n$. Thus $x_0$ is an $\varepsilon$-minimizer of $f$ on $\mathbb{R}^n$  if and only if $0\in \partial_{\varepsilon}f(x_0)$. For complete description of properties of $\varepsilon$-subdifferential see \cite{dhara2011optimality}.\\
The subdifferential defined above is only defined for convex functions, so the obvious question is to ask what about subdifferential of non-convex functions? We now discuss subdifferential of a non-convex function which is locally Lipschitz in nature. The relation of subdifferential and directional derivative as above becomes a key to develop the notion of a subdifferential for a locally Lipschitz functions.\\
 A function  $f:\mathbb{R}^n\rightarrow \mathbb{R}$ is Lipschitz around $x\in \mathbb{R}^n$, if there exists a neighborhood  $U_x$ of $x$ and $L_x\geq 0$ such that $\|f(y)-f(z)\|\leq L_x\|y-z\|,$ for all $y,z\in U_x$. The constant $L_x$ is the Lipschtiz constant of the function $f$ at the point $x$. A function $f$ is said to be locally Lipschitz if $f$ is Lipschitz around $x$ for any $x\in \mathbb{R}^n$. We shall focus in this article on MOP with locally Lipschitz objective and constraint functions. We now define the Clarke directional derivative of locally Lipschitz function $f$ at $x$ and in the direction $h\in \mathbb{R}^n$ as 
$$f^{\circ}(x,h)=\limsup_{y\to x,t\downarrow 0} \frac{f(y+t h)-f(y)}{t}.$$
The Clarke subdifferential of $f$ at $x\in \mathbb{R}^n$ is given as, 
$$\partial^{\circ} f(x)=\{\xi\in \mathbb{R}^n:~f^{\circ}(x,h)\geq \langle \xi, h \rangle,~\text{for all}~h\in \mathbb{R}^n\}.$$
For each $x\in \mathbb{R}^n,$ the set $\partial^{\circ} f(x)$ is non-empty, convex and compact. It is important to note that when function $f$ is convex, then $\partial^{\circ} f(x)=\partial f(x),$ for all $x\in \mathbb{R}^n$. Same as subdifferential for convex function, Clarke subdifferential has lots of nice properties. If $x_0\in \mathbb{R}^n$ is a local minimum of $f$ over $\mathbb{R}^n$, then $0\in \partial^{\circ}f(x_0)$ (for proof see \cite{rockafellar2015convex}). It also satisfy sum rule but it gives only one side containment, \textit{i.e.,} for given two locally Lipschtiz function $f$ and $g$, we have
$\partial^{\circ} (f+g)(x)\subset \partial^{\circ} f(x)+\partial^{\circ} g(x).$
\section{Approximate KKT conditions}

In this section, we begin by defining a notion of  modified $\epsilon$-KKT points which suits very well for the purpose of convex vector optimization problem. This notion is motivated by a similar notion defined in \cite{dutta2013approximate} for scalar optimization problem and \cite{durea2011stability} for convex vector optimization.
\begin{defn}\label{d21}
A feasible point $x_0\in X$ is said to be a modified $\varepsilon$-KKT point of MOP if for a given $\varepsilon \in \mathbb{R}_+$,  there exists $x_{\varepsilon}$ such that $\|x_0-x_{\varepsilon}\|\leq \sqrt{\varepsilon}$ and there exists $u_i\in \partial^{\circ} f_i(x_{\varepsilon})$ for all $i\in I$, $v_r\in \partial^{\circ} g_r(x_{\varepsilon})$ for all $r\in L$, vectors $\lambda \in \mathbb{R}_+^m$ with $\|\lambda\|=1$ and $\mu \in \mathbb{R}_+^l$ such that
\begin{eqnarray*}
\left\|\sum\limits_{i\in 1}^m \lambda_i u_i+\sum\limits_{r=1}^l \mu_r v_r \right\|\leq \sqrt{\varepsilon}, \text{ and, }\sum\limits_{r=1}^r\mu_r g_r(x_0)\geq -\varepsilon. \label{eq320}
\end{eqnarray*}
\end{defn}
\noindent Now we are in position to answer the first two questions asked in the introduction. To begin with we state two constraint qualifications, Slater constraint qualification (SCQ for short) and Basic constraint qualification (BCQ for short) which are used in the main results of this article (see \cite{rockafellar2009variational}).  
\begin{defn}\label{d22}
The MOP with constraint functions $g_r$ for all $r\in L$ to be convex satisfies Slater constraint qualification if there exists $\hat{x}\in X$ such that $g_r(\hat{x})<0$, for all $r\in L$. 
\end{defn}
\begin{defn}
The MOP with locally Lipschitz constraint functions $g_r$ for all $r\in L$ satisfies Basic Constraint Qualification (BCQ) at a point $\bar x$ if there exists no $p\in \mathbb{R}^l_+\setminus \{0\}$ such that $0\in  \sum\limits_{r\in L} p_r\partial^{\circ} g_r(\bar x)$.
\end{defn}
The next theorem answers the first question (\textbf{Q1}) raised in this article which says that if a sequence $x_k$ of modified $\varepsilon_k$-KKT points of MOP converges to a point $x_0$ where basic constraint qualification holds at $x_0$, then $x_0$ is a KKT point of MOP. Observe that we do not need convexity of the objective functions to prove the following result whereas we only require Lipschitz continuity of the objectives. It is important to note that the similar kind of result has been discussed under convexity assumption of objective functions in \cite{durea2011stability}.

\begin{thm}\label{t21}
Consider the problem MOP with locally Lipschitz data and let $\{\varepsilon_k\}$ to be a decreasing sequence of positive real numbers such that $\varepsilon_k \rightarrow 0$ as $k\rightarrow \infty$. Consider $\{x^k\}$ to be a sequence of feasible points of MOP with $x^k \rightarrow x_0$ as $k\rightarrow \infty$. Assume that for each $k$, $x^k$ is a modified $\varepsilon_k$-KKT point of MOP. Further, assume that the BCQ holds at $x_0$. Then $x_0$ is a KKT point of MOP.
\end{thm}
\textit{Proof:}
Note that $x^k$'s are feasible points, \textit{i.e.,} $g_r(x^k)\leq 0$ for all $r\in L$ and $x^k\rightarrow x_0$. Hence, using the convexity of $g_r$'s, we conclude that $g_r(x_0)\leq 0$ for all $r\in L$. Hence, $x_0$ is a feasible point of MOP. Now, as $x^k$ is a modified $\varepsilon_k$- KKT point, for each $k$,  Definition~\ref{d21} gives the existence of a point  $\hat{x^k}$ such that $\|x^k-\hat{x^k}\|\leq \sqrt{\varepsilon_k}$, the existence of  $u_i^k\in \partial^{\circ} f_i(\hat{x^k})$ and $v_r^k\in \partial^{\circ} g_r(\hat{x^k})$ for all $i \in I$ and $r\in L$, and the vectors $\lambda^k \in \mathbb{R}_+^m$  and $\mu^k \in \mathbb{R}_+^l$ with $\|\lambda^k\|=1$ such that  
\begin{eqnarray}
&& \left\|\sum\limits_{i\in I} \lambda^k_i u^k_i+\sum\limits_{r\in L} \mu^k_rv_r^k\right\|\leq \sqrt{\varepsilon_k},\label{eq01} \text{ and } \\
 && \sum\limits_{r\in L} \mu_r^k g_r(x^k)\geq -\varepsilon_k. \label{eq02}
\end{eqnarray}
We first claim that $\{\mu^k\}$ is bounded. To prove our claim,  on the contrary assume that $\{\mu^k\}$ is unbounded. Thus,  $\|\mu^k\|\rightarrow \infty$ as $k\rightarrow \infty$. Further, Equation~(\ref{eq01}), can be re-written as 
\begin{equation}\label{eq03}
\left\|\sum\limits_{i\in I} \frac{\lambda^k_i}{\|\mu^k\|} u^k_i+\sum\limits_{r \in L} \frac{\mu^k_r}{\|\mu^k\|} v_r^k \right\|\leq \frac{1}{\|\mu^k\|}\sqrt{\varepsilon_k}.
\end{equation}
Then, in Equation~\eqref{eq03}, we observe the following:
\begin{enumerate}
\item As $\varepsilon_k$ converges to $0$, the same holds for $\frac{1}{\|\mu^k\|}\sqrt{\varepsilon_k}$.
\item Let $p_r^k=\frac{\mu_r^k}{\|\mu^k\|} \in \mathbb{R}_+$, for all $r\in L$. As $\|p^k\| =1$,  $\{p^k\}$ is a bounded sequence. So, by the Bolzano-Weierstrass theorem, there exists a subsequence of $\{p^k\}$ which converges to $\hat{p}\in \mathbb{R}^l_+$  with $\|\hat{p}\| = 1$.  In fact, without loss of generality, we can assume that  $p^k_r$ converges to $\hat{p}_r$. Hence, for all $r\in L$
\begin{equation}\label{eq04}
\frac{\mu^k_r}{\|\mu^k\|}=p^k_r\rightarrow \hat{p}_r,~as ~k\rightarrow \infty. 
\end{equation}
\item As  $f_i$'s are locally Lipschitz functions,  their Clarke subdifferential  are locally bounded, \textit{i.e.,} for $x_0\in X$, there exists $\delta>0$ such that for all $z \in B_{\delta}(x_0)$, $\partial^\circ f_i(z)\subset K_i$, where,  for all $i\in I$, $K_i$'s are bounded sets on $\mathbb{R}^n$. Since $x^k\rightarrow x_0$, there exists $k_0\in \mathbb{N}$ such that, for all $k\ge k_0 $, $x^k\in B_{\delta}(x_0)$. Therefore, by choosing $K=\bigcup\limits_{i\in I}{\tilde{K_i}}$ where $\tilde{K_i}=  K_i\cup\partial^{\circ} f_i({x^1})\cup \ldots \partial^{\circ} f_i({x^k_0})  $, we get $\partial^{\circ} f_i({x^k})\subset K$, for all $i\in I$ and $k\ge 0$. Hence, the sequence $\{u_i^k\}$, where $u_i^k\in \partial^{\circ} f_i({x^k})$, is bounded for all $i \in I$. Hence, using the fact that $\|\lambda^k\|=1$ and $\|\mu^k\|\rightarrow \infty$, we deduce that for all $i\in I$,
\begin{eqnarray}\label{eq06}
\frac{\lambda_i^k}{\|\mu^k\|}u_i^k\rightarrow 0, ~as ~k\rightarrow \infty.
\end{eqnarray}

\item An argument similar to the previous part implies that the sequence $\{ v_r^k\}$ where $v_r^k\in \partial^{\circ} g_r(\hat{x^k}) $, for each fixed $r \in L$, is bounded. Hence, 
 the sequence $\{v_r^k\}$ has a limit point, for all $r\in L$, say $\hat{v}_r$. Without loss of generality, we can assume that for all $r\in L$, 
\begin{equation}\label{eq07}
v_r^k\rightarrow \hat{v}_r, \text{ as } k\rightarrow \infty.
\end{equation}
Since $\partial^{\circ} g_r$'s is graph closed and $\hat{x}_k\rightarrow x_0$, one has $\hat{v}_r\in \partial^{\circ} g_r(x_0)$ for all $r\in L$.
\end{enumerate}
Now, take the limit as $k\rightarrow \infty $ in Inequality~(\ref{eq03}) and  in view of the above observations~(\ref{eq04}),(\ref{eq06}) and (\ref{eq07}), we get, 
$$\left\|\sum_{r\in L} \hat{p}_r\hat{v}_r\right\| \leq 0.$$  
Hence, we have $\sum\limits_{r\in L} \hat{p}_r \hat{v}_r=0$, where $\hat{p}\in \mathbb{R}^l_+$ with $\|\hat{p}\|=1$ and $\hat{v}_r\in \partial^{\circ} g_r(x_0)$ for all $r\in L$. This contradicts the assumption that BCQ holds at $x_0$. Therefore, we have shown the correctness of our claim, {\textit i.e.}, the sequence $\{\mu^k\}$ is a bounded.

As $\{\mu^k\}$ is a bounded sequence, an argument similar to the one above, implies that there exist $\hat{\mu}\in \mathbb{R}^l_+$ such that $\mu_k \rightarrow \hat{\mu} ~as ~k\rightarrow \infty$. Similarly,  the sequences $\{\lambda^k\}$ and $\{u^k\}$ have  limit points, say $\hat{\lambda}$ and $\hat{u}$, respectively, with $\|\hat{\lambda}\|=1$ and  $\lambda_k\rightarrow \hat{\lambda}$, $u^k \rightarrow \hat{u}$.
Now taking $k\rightarrow \infty $ in Inequality~(\ref{eq01}), we get
$\|\sum\limits_{i \in I} \hat{\lambda_i} \hat{u}_i+\sum\limits_{r\in L} \hat{\mu_r} \hat{v}_r\|\leq 0.$
Thus  
\begin{eqnarray}\label{eq09}
\sum\limits_{i \in I} \hat{\lambda_i} \hat{u}_i+\sum\limits_{r \in L} \hat{\mu_r} \hat{v}_r= 0, &&  \hspace*{-.2in}
\text{ where } \hat{\lambda}\in \mathbb{R}^m_+ \text{ with } \|\hat{\lambda}\|=1, \nonumber \\ 
&&  \hat{\mu}\in \mathbb{R}^l_+, \hat{u}_i\in \partial^{\circ} f_i(x_0) \text{ and } \hat{v}_r\in \partial^{\circ} g_r(x_0).
\end{eqnarray}
Since, $x_0$ is a feasible point of MOP and $\hat{\mu}_r\geq 0$ for all $r\in L$, we have
$\sum\limits_{r \in L} \hat{\mu}_r g_r(x_0) \leq 0.$
Taking $k\rightarrow \infty$ in Inequality~(\ref{eq02}), we get $\sum\limits_{r\in L} \hat{\mu}_r g_r(x_0) \geq 0$ and thus, we conclude that 
\begin{equation}\label{eq012}
\sum\limits_{r\in L} \hat{\mu}_rg_r(x_0) =0.
\end{equation}
The Inequalities~(\ref{eq09}) and (\ref{eq012}) together imply that $x_0$ is a KKT point of MOP.
\hfill$\Box$\\

\noindent The next theorem deals with the second question asked in the article. Basically, \textbf{Q2} for multiobjective problem can be framed as follows: for every local Pareto points of MOP, does there exists sequence which converges to the point and the sequence has a subsequence which satisfies some type of approximate KKT conditions? We answer this question in Theorem \ref{t22} for MOP with locally Lipschitz objective function and Slater constraint qualification. This result shows that we always have a sequence converging to a local Pareto point of MOP with approximate KKT type of conditions which implies that the idea of constructing approximate KKT type conditions is essential in multiobjective theory. Note that \textbf{Q2} has not been addressed in \cite{durea2011stability}. Before we state theorem, we present the following lemma which will be needed in the proof of the result. This lemma is a special case of Theorem 2.44 in \cite{mordukhovich2013easy}. 
 \begin{lem}\label{lemma}
 Let $A$ and $B$ be two non-empty subsets of $\mathbb{R}^m$. Let $A\cap B\ne \emptyset$ and $\bar{x}\in A\cap B$. Assume that the following qualification condition holds:
 $$N_A(\bar{x})\cap (-N_B(\bar{x}))=\{0\}.$$
 Then, $N_{A\cap B}(\bar{x})=N_A(\bar{x})+N_B(\bar{x})$. 
 \end{lem}

\begin{thm}\label{t22}
Consider the problem MOP with locally Lipschtiz objectives $f_i$'s for all $i\in I$ and $g_r$'s for all $r \in L$ to be a convex functions which satisfies the Slater constraint qualification. Further, assume that $x_0$ is a local weak Pareto minima and consider $\{\varepsilon_k\}$ to be a decreasing sequence of positive real numbers converging to $0$. Then, there exists a sequence $\{x^k\}$ of feasible points converging to $x_0$ which has a subsequence $\{y^k\}$ of $\{x^k\}$ such that for each $y^k$, there exists $\hat{y}^k$ satisfying
\begin{enumerate}
\item  $\|y^k-\hat{y}^k\|\leq \sqrt{\varepsilon_k}$, 
\item there exists $u_i^k\in \partial^{\circ} f_i(\hat{y}^k)$ and $v_r^k\in \partial^{\circ} g_r(\hat{y}^k)$, for all $i \in I$ and $r\in L$, such that 
\begin{eqnarray}
\left\|\sum\limits_{i\in I} \lambda^k_i u^k_i+\sum\limits_{r\in L} \mu^k_r v_r^k \right\|\leq \sqrt{\varepsilon_k}, \label{eq311}\\
\sum\limits_{r\in L}\mu_r^k g_r(\hat{y}^k)= 0,\label{eq312}
\end{eqnarray}
where  $\lambda^k \in \mathbb{R}^m_+$ with $\|\lambda^k\|=1$ and $\mu^k\in \mathbb{R}^l_+$. 
\end{enumerate}

\end{thm}
\textit{Proof:} By  assumption, $x_0$ is a locally Pareto minimizer of MOP, \textit{i.e.}, there exists $\delta>0$ such that 
\begin{equation*}
f(x)-f(x_0)\not \in - \mathbb{R}^m_+\setminus \{0\}, \text{ for all } x\in V,
\end{equation*}
equivalently, 
\begin{equation}\label{eq013}
f(x)-f(x_0)\in \tilde{W}, \text{ for all } x\in V,
\end{equation}
where $\tilde{W}:=\mathbb{R}^m\setminus (- \mathbb{R}^m_+\setminus \{0\})$ and $V=X\cap \overline {B_{\delta}(x_0)}$. The convexity of the constraint functions $g_r$'s together with closed convex feasible set $X$ implies that  $V$ is a closed, convex and bounded set. As $x_0\in V$, there exists a sequence $x^k$ in $X$ with $x^k$ converging to $x_0 \in V$ and $x^k \in V$, 
for all $k$ sufficiently large. We have broken the rest of the proof in two steps. For the first step, we prove that there exists a sub-sequence $\{y^k\}$ of $\{x^k\}$ such that  $y^k \in V$ and is an $\varepsilon_k e$-Pareto minima of MOP with feasible set as $V$ where $e=(1,\dots,1)^T$ and $\varepsilon_k>0$. 

As  $f_i$'s, for $i\in I$, are locally Lipschitz, $f_i(x^k) \rightarrow f_i(x_0)$ as $k\rightarrow \infty$, for all $i\in I$. So, for a given $\varepsilon_1 >0$,  for each $i \in I$ there exist  natural numbers $N^i_1$, such that 
$$| f_i(x^k)-f_i(x_0) |< \varepsilon_1, \text{ for all } k \ge N^i_1.$$
Now choose $N_1=\max\{N^1_1,N^2_1,\ldots,N^m_1\}$. Thus, for all  $i\in I$ 
\begin{equation}\label{eqn:epsi1}
| f_i(x^k)-f_i(x_0)|< \varepsilon_1, \text{ for all } k \ge N_1.
\end{equation}
Choose $y^1=x^{N_1}$, then $| f_i(y^1)-f_i(x_0) |< \varepsilon_1,$ or equivalently, 
\begin{equation}\label{eq015}
f(x_0)+e\varepsilon_1-f(y^1)\in {int}(\mathbb{R}^m_+).
\end{equation}
Note that $\tilde{W} +int(\mathbb{R}^m_+)\subseteq \tilde{W}$, hence, (\ref{eq013}) and~(\ref{eq015}) together gives 
\begin{equation}
f(x)+e\varepsilon_1-f(y^1) \not \in  - \mathbb{R}^m_+\setminus \{0\}, \text{ for all } x\in V. 
\end{equation} 
Take $\varepsilon_2 < \varepsilon_1 $ and a similar argument applied to the sequence $\{x_{N_1}, x_{N_1 + 1}, x_{N_1 + 2}, \ldots\}$ gives an element $y^2=x^{N_2}$, with $N_2 > N_1$, such that 
$ f(x)+\varepsilon_2 e-f(y^2)\not \in  - \mathbb{R}^m_+\setminus \{0\}\text{ for all }x\in V.$ Proceeding as above, gives a sub-sequence $\{y^k\}$ of $\{x^k\}$ such that  $y^k \in V$ and 
\begin{equation}\label{eq016}
f(x)+\varepsilon_ke-f(y^k)\not \in - \mathbb{R}^m_+\setminus \{0\}, \text{ for all }x\in V.
\end{equation}
Hence, $y^k \in V$ is an $\varepsilon_k e$-Pareto minima of MOP with feasible set as $V$. This  completes the proof of the first step. We now come to the second step to complete the proof. 

 Since each  $f_i$ is locally Lipschitz,  $f$ is ($e,\mathbb{R}^m_+$)-lower semi continuous and $\mathbb{R}_+^m$-bounded below. Thus,  the vector Ekeland Variational Principle (Theorem \ref{d17}) gives the existence of $\hat{y}^k \in V$, for each $y^k \in V$, such that $\|\hat{y}^k-y^k\|\leq \sqrt{\varepsilon_k},$ and for all $x\in V\setminus \{\hat{y}^k\},$
\begin{enumerate} 
\item  $f(x)+\varepsilon_k e-f(\hat{y}^k)\not \in -int (\mathbb{R}_+^m)$, and
\item   $f(x)+\sqrt{\varepsilon_k}\|\hat{y}^k- y^k\| e-f(\hat{y}^k)\not \in -int (\mathbb{R}^m_+)$.
\end{enumerate}
Thus from above, we conclude that $\hat{y}^k$ is a weak Pareto minimizer of the problem 
$$\min\limits_{x\in V} g(x), \text{ where } g(x) = f(x)+\sqrt{\varepsilon_k}\|x-\hat{y}^k\| e.$$
Now, using the necessary optimality condition for the above multiobjective problem, there exists $\lambda^k\in \mathbb{R}^m_+$ with $\|\lambda^k\|=1$ such that
$$0\in \sum\limits_{i\in I}\lambda_i^k \partial^{\circ}g_i(\hat{y}^k)+N_V(\hat{y}^k),$$
where $N_V(\hat{y}^k)$ is the normal cone to the set $V$ at $\hat{y}^k$. For proof of above result see for example, page 137 of Chapter 5 in \cite{dutta2012strong}.
 Now applying sum rule for the Clarke subdifferential (see \cite{clarke1990optimization}) and using the fact that subdifferential of the norm function at origin is the unit ball, we get
\begin{equation}\label{eq018}
0\in \sum\limits_{i\in I}\lambda_i^k \partial^{\circ} f_i(\hat{y}^k)+\sqrt{\varepsilon_k}B_{1}(0)+N_V(\hat{y}^k).
\end{equation}
Since $x^k\rightarrow x_0$ and $y^k$ is a sub-sequence of  $\{x^k\}$,   $y^k\in X\cap B_{\delta}(x_0)$, for  sufficiently large $k$. As  $\hat{y}^k\in B_{\sqrt{\epsilon_k}}(y^k)$ and $\epsilon_k \rightarrow 0$, for sufficiently large $k$,  $B_{\sqrt{\epsilon_k}}(y^k) \subset B_{\delta}(x_0)$.  Hence, $\hat{y}^k\in B_{\delta}(x_0)$, for  $k$ sufficiently large.

Clearly,  $X\cap \overline{B_{\delta}(x_0)}\not = \emptyset$. We will now see that the qualification condition for Lemma \ref{lemma} holds in this case. Since $\hat{y}^k\in B_{\delta}(x_0)$, we see that   $\hat{y}^k\in{int} \overline{B_{\delta}(x_0)}$, thus $N_{\overline{B_{\delta}(x_0)}}(\hat{y}^k)=\{0\}$. Hence, $N_V(\hat{y}^k)\cap(-N_{\overline{B_{\delta}(x_0)}}(\hat{y}^k))=\{0\}$. Therefore, using Lemma \ref{lemma}, we conclude that 
$$ N_V(\hat{y}^k)=N_{X\cap \overline{B_{\delta}(x_0)}}(\hat{y}^k)
=N_X(\hat{y}^k)+N_{\overline{B_{\delta}(x_0)}}(\hat{y}^k).$$
Thus $N_V(\hat{y}^k)=N_X(\hat{y}^k).$ Hence, we can rewrite (\ref{eq018}) as 
\begin{equation}\label{eq019}
0\in \sum\limits_{i\in I}\lambda_i^k \partial^{\circ} f_i(\hat{y}^k)+\sqrt{\varepsilon_k}B_{1}(0)+N_X(\hat{y}^k).
\end{equation} 
Further as the Slater constraint qualification holds, using Corollary~$23.7.1$ of \cite{rockafellar2015convex}, 
$$
N_X(\hat{y}^k)=\{\sum\limits_{r\in L}\mu_r^kv_r^k: v_r\in \partial g_r(\hat{y}^k), ~\mu^k_r\geq 0, ~\mu_r^kg_r^k(\hat{y}^k)=0,~r\in L\}.
$$
Now using the above form of $N_X(\hat{y}^k)$ and (\ref{eq019}), it is evident that there exists $u_i^k\in \partial^{\circ} f_i(\hat{y}^k)$ for all $i\in I$, $v_r^k\in \partial g_r(\hat{y}^k)$ for all $r\in L$ and scalars $\lambda^k \in \mathbb{R}^m_+$ with $\|\lambda^k\|=1$, $\mu^k\in \mathbb{R}^l_+$ such that \eqref{eq311} and (\ref{eq312}) holds. This completes the proof of the second part and hence  the proof of the theorem is complete.\hfill$\Box$

\begin{remark}\label{rm2}
In the above theorem, the objective functions are taken to be locally Lipschitz only. If the objective function $f_i$'s are convex as well, then we have a more concrete result. To proof the next result we need the following Lemma \ref{lem01} and a result from \cite{durea2011stability} which will play a key role in proving the Theorem \ref{c3t6}.
\end{remark}

\begin{lem}\label{lem01}
Consider the problem MOP with each objective functions $f_i$'s and constraint function $g_r$'s to be convex. Then every local Pareto minima is a global Pareto minima. 
\end{lem}
\begin{thm}[Theorem 3.6 of \cite{durea2011stability}]\label{c3co6}
Let $x_0$ be a $\varepsilon e$-weak Pareto minima of the problem MOP with each $f_i$'s and $g_r$'s to be convex functions and assume that Slater constraint qualification holds. Then $x_0$ is a modified $\sigma$-KKT point where $\sigma\in (0,\|e\|\varepsilon]$. 
\end{thm}
\begin{thm}\label{c3t6}
Consider the problem MOP with each $f_i$ and $g_r$  being convex functions, for all $i\in I$ and $r\in L$. Let $x_0$ be a Pareto minima and let the Slater constraint qualification hold. Then, for decreasing sequence of positive real numbers  $\{\varepsilon_k\}$ converging to $0$, there exists a feasible sequence $\{x^k\}$ converging to $ x_0$  and a sub-sequence $\{y^k\}$ of $\{x^k\}$ such that each $y^k$ is a modified $\sigma_k$-KKT point with $\sigma_k\in (0, \; \|e\|\varepsilon_k ]$.
\end{thm}
\textit{Proof:}
Since the problem data is convex, local Pareto point is global. Now proceed as in  the proof of Theorem~\ref{t22} to get a sub-sequence $\{y^k\}$ of $\{x^k\}$ such that $y^k$  is a $\varepsilon_k e$-Pareto minima of MOP with feasible set as $V$), where $V=X\cap \overline{B_{\delta}(x_0)}$ with $\delta >0$, \textit{i.e.,} $y^k$ is a local $\varepsilon_k e$-Pareto minima of MOP. So, by using the assumption of convexity and Lemma~\ref{lem01}, we conclude that $y^k$ is a $\varepsilon_k e$-Pareto minima of MOP. Now using Theorem~\ref{c3co6}, we conclude that $y^k$ is a modified $\sigma_k$-KKT point with $\sigma_k\in (0, \; \| e\|\varepsilon_k]$.\hfill$\Box$

\section{Approximate $\hat{M}$-Geoffrion solutions, Saddle points, and KKT conditions}
In this section, we analyze saddle point conditions and KKT type conditions for the $(\hat{M},\epsilon)$-Geoffrion solutions which give a complete characterization of the considered proper points. We also discuss a scalarization rule for the $(\hat{M},\epsilon)$-Geoffrion solutions which is a connecting bridge for deducing saddle point and KKT type conditions. Before discussing the mentioned results, we shall observe that there is a characterization of $(\hat{M},\epsilon)$-Geoffrion proper points by the system of inequalities which appeared in \cite{shukla2019practical}. For a given $\epsilon\in \mathbb{R}^m_+$ and $\hat{M}>0$, consider $x_0\in X$, $i\in I$ and define the following system of inequalities ($\mathcal{Q}_i(x_0)$) as

\begin{eqnarray*}
\left\{ \begin{array}{ll}
       -f_i(x_0)+ f_i(x)+\epsilon_i<0,\\
-f_i(x_0)+ f_i(x)+\epsilon_i<\hat{M} (f_j(x_0)- f_j(x)-\epsilon_j),\text{ for all } j\in I\setminus\{i\}\\
x\in {X}.\end{array} \right.
\end{eqnarray*}

\begin{prop}\label{pro1}
For given $\epsilon\in \mathbb{R}^m_+$ and  $\hat{M}>0$, consider the problem MOP. Then a point $x_0\in \mathcal{G}_{\hat{M},\epsilon}(f,X)$ if and only if for each $i\in I$, the system $\mathcal{Q}_i(x_0)$ is inconsistent. 
\end{prop}

\noindent The above proposition follows from the definition of Proof of the $(\hat{M},\epsilon)$-Geoffrion proper solutions, for complete proof, see \cite{shukla2019practical}. Before discussing the saddle point conditions for the $(\hat{M},\epsilon)$-Geoffrion proper solutions, let us discuss the correspondence between $(M,\epsilon)$-Geoffrion proper solutions and solution of the weighted sum scalar problem. As mentioned earlier, this correspondence plays a pivotal role to prove main results of this section. To this end, let for $s^\ast\in\mathbb{R}^m_+$, the weighted sum scalar problem $P(s^*)$ be defined as $\min\limits_{x\in X} ~\langle s^*,f(x)\rangle.$
\begin{thm}\label{t2}
For a given $\epsilon \in \mathbb{R}^m_+$, $\hat{M}>0$, let $x_0$ is a $\langle s^\ast,\epsilon\rangle$-minimum of $P(s^\ast)$, where $s^\ast \in \text{int}( \mathbb{R}^m_+)$. If $\hat{M}\geq (m-1)\max\limits_{i,j}\{\frac{s^\ast_i}{s^\ast_j}\}$, then $x_0$ is a $(\hat{M},\epsilon)$-Geoffrion proper solution of MOP, \textit{i.e.,} $x_0\in \mathcal{G}_{\hat{M},\epsilon}(f,X)$.

\end{thm}
\textit{Proof:}
Let us assume on the contrary that $x_0\notin \mathcal{G}_{\hat{M},\epsilon}(f,X)$. Therefore, from Proposition~\ref{pro1} we obtain an ${i}\in I$ such that  $\mathcal{Q}_i(x_0)$ is consistent. Without loss of generality, we assume that $i=1$. Thus, the system $\mathcal{Q}_i(x_0)$, written as
\begin{eqnarray*}\label{eqn}
\left\{ \begin{array}{ll}
       -f_1(x_0)+ f_1(x)+\epsilon_1<0,\\
-f_1(x_0)+ f_1(x)+\epsilon_1<\hat{M} (f_j(x_0)- f_j(x)-\epsilon_j), \quad j\in I\setminus\{1\}\\
x\in {X}.\end{array} \right.
\end{eqnarray*}
has a solution.  As $\hat{M}\ge (m-1)\{\frac{s^\ast_j}{s^\ast_i}\}$ for all $s^\ast\in  \text{int}(\mathbb{R}^m_+)$, the consistency of system $\mathcal{Q}_i(x_0)$ implies that
\begin{equation*}\label{eqnn}
  s^\ast_1(-f_1(x_0)+ f_1(x)+\epsilon_1)<s^\ast_j (m-1) (f_j(x_0)- f_j(x)-\epsilon_j),\text{ for all } j\in I\setminus\{1\}.
\end{equation*}
Summing the above equation for all $j\in I\setminus\{1\}$, we obtain that
\begin{equation*}
  s^\ast_1(-f_1(x_0)+ f_1(x)+\epsilon_1)<\sum_{j=2}^ms^\ast_j  (f_j(x_0)- f_j(x)-\epsilon_j),
\end{equation*}
which further implies
\begin{equation}\label{eqnnn}
\langle s^*,f(x_0)\rangle - \langle s^*, f(x)\rangle -\langle s^*,\epsilon \rangle>0.
\end{equation}
Since (\ref{eqnnn}) is a contradiction to the $\langle s^\ast,\epsilon\rangle$-minimality of $P(s^\ast)$. Therefore, the theorem follows.\hfill$\Box$\\

\noindent All the solutions from $\mathcal{G}_{\hat{M},\epsilon}(f,X)$ satisfy an upper trade-off bound of $\hat{M}$ (in the sense of Geoffrion-proper efficiency). Smaller bounds are more relevant to the decision maker as they provide tighter trade-offs among the criteria values. Therefore, it is of interest to find the minimum ${M}$ such that $\mathcal{G}_{{M},\epsilon}(f,X)$ is non-empty.  Under the conditions of Theorem~\ref{t2}, we need minimum value of $\hat{M}$ equals $m-1$, and this occurs when all components of $s^\ast$ are identical. The next example shows that if conditions in Theorem~\ref{t2} are not satisfied, then even smaller values of $\hat{M}$ are possible. This is the case with non-convex or discrete multicriteria optimization problems. In the following example, we consider $\epsilon=0$ and find $\hat{M}$-Geoffrion proper points.
\begin{example}\rm
Let $X:=\{(0,0,1)^\top,\,(0,1,0)^\top,\,(1,0,0)^\top,(1/\sqrt{3},1/\sqrt{3},1/\sqrt{3})^\top\}$, $m=3$, and $f$ be the identity mapping. The sets $\mathcal{G}_{2}(f,X)$ and $\mathcal{G}_{1}(f,X)$ can be easily computed as follows:
\begin{eqnarray*}
\mathcal{G}_{2}(f,X)&=&\{(0,0,1)^\top,\,(0,1,0)^\top,\,(1,0,0)^\top,(1/\sqrt{3},1/\sqrt{3},1/\sqrt{3})^\top\},\\
\mathcal{G}_{1}(f,X)&=&\{(0,0,1)^\top,\,(0,1,0)^\top,\,(1,0,0)^\top\}.
\end{eqnarray*}
Moreover, $\mathcal{G}_{M}(f,X)=\emptyset$ for $M<1$. Therefore, the minimum value of ${M}$ is 1.
\end{example}

\noindent The converse of Theorem \ref{t2} also holds with convexity assumption on the objective functions and the feasible set. Since, if for each $r\in L$, $g_r$ is convex, then the feasible set $X$ is a convex set. We have the following result.
\begin{thm}\label{t3}
Let us consider the problem MOP where for each $i\in I$ and $r\in L$, $f_i$ and $g_r$ are convex functions. If $x_0\in \mathcal{G}_{\hat{M},\epsilon}(f,X)$, then there exists an $s^*\in {\rm int}(\mathbb{R}^m_+)$ such that $x_0$ is a $\langle s^*,\epsilon\rangle$-minimum of $P(s^*)$.
\end{thm}
\textit{Proof:}
Let $x_0\in \mathcal{G}_{\hat{M},\epsilon}(f,X)$. Then using Proposition \ref{pro1}, we obtain that the system $\mathcal{Q}_i(x_0)$  is inconsistent, for each $i\in I$. Applying the Gordan's Theorem of the alternative (see \cite{rockafellar2015convex}), we conclude, after some rearrangements, that for each $i\in I$, there exists scalars 
$\lambda_j^i\geq 0$ with $\sum\limits_{j \in I} \lambda_j^i=1$ such that, for all $x\in X$
\begin{eqnarray*}
f_i(x)+\hat{M}\sum_{j\in I,j\not = i} \lambda_j^if_j(x)\geq f_i(x_0)+\hat{M}\sum_{j\in I,j\not = i}\lambda_j^i f_j(x_0)-\left[\epsilon_i+\hat{M}\sum_{j\in I,j\not = i}\lambda_j^i\epsilon_j \right].
\end{eqnarray*}
 Therefore, by summing over all $i$, we get 
\begin{eqnarray*}
\sum_{i\in I} f_i(x)+\hat{M}\sum_{i\in I}\sum_{j\in I,j\not = i}\lambda_j^i f_j(x) &\geq& \sum_{i\in I} f_i(x_0)+\hat{M}\sum_{i\in I} \sum_{j\in I,j\not = i}\lambda_j^if_j(x_0)\\
&& \hspace{.35in} -\sum_{i\in I}\left[\epsilon_i+\hat{M}\sum_{j\in I,j\not = i}\lambda_j^i\epsilon_j \right].
\end{eqnarray*}
Hence, for all $x\in X$,
{\small \begin{eqnarray*}
\sum_{j\in I}\left[1+\hat{M}\sum_{i\in I,i\not = j}\lambda_j^i\right]f_j(x) \geq \sum_{j\in I}\left[1+\hat{M}\sum_{i\in I,i\not = j} \lambda_j^i\right]f_j(x_0)-\sum_{j\in I}\left[1+\hat{M}\sum_{i\in I,i\not = j}\lambda_j^i\right]\epsilon_j.
\end{eqnarray*}}
Setting $s_j=1+\hat{M}\sum\limits_{i\in I,i\not = j} \lambda_j^i$, gives  $s\in { int}(\mathbb{R}^m_+)$ and  $x_0$ is a $\langle s,\epsilon\rangle$-minimum of $P(s)$.\hfill$\Box$
\begin{remark}\rm
 Theorem \ref{t3} can also be proved by noting the fact that each $(\hat{M},\epsilon)$-Geoffrion proper point is $\epsilon$-Geoffrion proper point with constant $\hat{M}>0$. Hence using Theorem 3.15 form \cite{EHR}, we can deduce the above result. Now if we denote the set of $\langle s^*,\epsilon\rangle$-minimum of $P(s^*)$ by $Sol_\epsilon(P(s^*))$, then Theorem \ref{t2} and \ref{t3} implies that under convexity assumption on data and for a given $\hat{M}$, there exists $s^*\in \text{int}(\mathbb{R}^m_+)$ such that
$$ Sol_{\epsilon}(P(s^*))\subseteq \mathcal{G}_{\hat{M},\epsilon}(f,X) \subseteq \bigcup\limits_{s\in \text{int} (\mathbb{R}^m_+)}Sol_{\epsilon}(P(s)).$$
\end{remark}

\noindent Now we come to the main attraction of this section, the saddle point conditions for $(\hat{M},\epsilon)$-Geoffrion proper solutions. For this study, we consider the problem MOP where each $f_i$, $i\in I$ and $g_j$, $j\in L$ are a convex function.  Whenever the data of problem is convex , we shall denote the problem MOP as CMOP. Given $\hat{M}>0$, and any index $i\in I$, we define the $(\hat{M},i)$-Lagrangian associated with CMOP as follows
\begin{equation}\label{eq410}
L^{\hat{M}}_i(x,\tau^i, \mu^i)=f_i(x)+\sum_{j\in I,j\not = i} \tau^i_j \hat{M} f_j(x)+ \sum_{r\in L}\mu_r^ig_r(x),
\end{equation}
where $\mu^i=(\mu^i_1,\mu_2^i,...,\mu^i_l) \in \mathbb{R}^l_+$ and $\tau^i=(\tau^i_1, \tau^i_2,...,\tau^i_m)\in S^m$ with $S^m=\{x\in \mathbb{R}^m: 0\leq x_i\leq 1,i\in I, \sum_{i=1}^mx_i=1\},$ the unit simplex in $\mathbb{R}^m$.
The motivation behind considering the above Lagrangian comes from the $i$th-objective Lagrangian problem defined in Chapter 4 of \cite{chankong2008multiobjective}. In \cite{chankong2008multiobjective}, they used the above Lagrangian form as a scalarization scheme of multiobjective problems. In the same spirit as \cite{chankong2008multiobjective}, we get a scalar structure of Lagrangian functions which is comparatively easy than vector-valued Lagrangian to work with. Our aim here is to show the key role played by the $(\hat{M},i)$-Lagrangian in analyzing and characterizing the Geoffrion ($\hat{M}, \epsilon$)-Proper solutions.

\begin{thm}\label{t041}
For a given $\epsilon \in \mathbb{R}^m_+$ and $\hat{M}>0$, let us consider the problem CMOP which satisfy the Slater constraint qualification. If $x_0\in \mathcal{G}_{\hat{M},\epsilon}(f,X)$ then for each $i$, there exists $\bar{\tau}^i\in S^m$, $\bar{\mu}^i\in \mathbb{R}^l_+$ such that for all $x\in \mathbb{R}^n$ and $\mu \in\mathbb{R}_+^m $,
\begin{description}
\item $(i)$ $L^{\hat{M}}_i(x_0,\bar \tau^i, \mu)-\bar \epsilon_i \leq L^{\hat{M}}_i(x_0,\bar \tau^i, \bar \mu^i)\leq L^{\hat{M}}_i(x,\bar \tau^i, \bar \mu^i)+\bar \epsilon_i$
\item($ii$) $\sum\limits_{r\in L} \bar \mu_r^i g_r(x_0) \geq -\bar \epsilon_i$, 
\end{description}
where $\bar \epsilon_i=\epsilon_i+\sum\limits_{j=1,\\ j\not = i}^m\tau^i_j \hat{M} \epsilon_j.$ Conversely if $x_0\in \mathbb{R}^n$ be such that for each $i\in I$, there exists $(\bar \tau^i, \bar \mu^i)\in S^m\times \mathbb{R}^l_+$ such that $(i)$ and ($ii$) holds then $x_0\in \mathcal{G}_{\tilde{M},2\epsilon}(f,X)$, where $\tilde{M}\geq(1+\hat{M})(m-1)$.
\end{thm}
\textit{Proof:}
It is evident from Proposition \ref{pro1} that if $x_0\in \mathcal{G}_{\hat{M},\epsilon}(f,X)$, then for each $i\in I$, the system $\mathcal{Q}_i(x_0)$, re-written as
\begin{align*}
 & -f_i(x_0)+ f_i(x)+\epsilon_i<0,\\
& -f_i(x_0)+ f_i(x)+\epsilon_i<M (f_j(x_0)- f_j(x)-\epsilon_j),\text{ for all } j\in I\setminus\{i\}\\
& g_r(x)\leq 0,~ r\in L 
\end{align*}
has no solution, for all $x\in \mathbb{R}^n$. It is easy to observe that the system $\mathcal{Q}_i(x_0)$ has no solution, if we replace $g_r\leq 0$ by $g_r<0$ for all $r\in L$. Now by applying the Gordan's theorem of the alternative (see \cite{rockafellar2015convex}), there exists $\tau^i=(\tau^i_1,\ldots,\tau^i_m)\in \mathbb{R}^m_+$ and $\mu^i=(\mu^i_1,\ldots,\mu^i_l)\in \mathbb{R}^l_+$ with $(\tau^i,\mu^i)\not = 0$ such that for all $x\in \mathbb{R}^n,$
\begin{multline*}
\tau^i_i(f_i(x)-f_i(x_0)+\epsilon_i)+
\sum\limits_{j\in I,j\not = i}\tau^i_j(f_i(x)+\hat M f_j(x) - f_i(x_0) \\
- \hat M f_j(x_0) +\epsilon_i
+\hat{M}\epsilon_j)
+\sum\limits_{r\in L} \mu^i_r g_r(x) \geq 0.
\end{multline*}
Hence, for all $x\in \mathbb{R}^n,$
\begin{eqnarray}\label{eq430}
\bigl(\sum\limits_{j\in I} \tau^i_j  \bigr)(f_i(x)-f_i(x_0)+\epsilon_i)
+\sum\limits_{j\in I,j\not = i} \left[\tau^i_j \hat M f_j(x)-  \tau^i_j \hat M f_j(x_0)
+\tau^i_j \hat M \epsilon_j \right]\nonumber \\ 
 + \sum\limits_{r\in L} \mu^i_r g_r(x) \geq 0.
\end{eqnarray}
Now, we first claim that $\tau^i=(\tau^i_1,\ldots,\tau^i_m)\not =0$. For if, $\tau^i=0$ then $\mu^i \ne 0$ and Inequality~\eqref{eq430} reduces to $\sum\limits_{r\in L} \mu^i_r g_r(x) \geq 0$, for all $x\in \mathbb{R}^n$.  But,   the Slater constraint qualification implies that there exists a point, say $\hat{x}\in \mathbb{R}^n$, such that 
$g_r(\hat{x})<0$. As $\mu^i \ne 0$ and $\mu^i \in \mathbb{R}^l_+$, we obtain
 $\sum\limits_{r\in L} \mu^i_r g_r(x) <0$, a contradiction to $\sum\limits_{r\in L} \mu^i_r g_r(x) \geq 0.$ Hence, $\tau^i \not = 0$ and thus $\sum\limits_{j\in I} \tau^i_j >0$. Thus, dividing  Inequality~(\ref{eq430}) by $\sum\limits_{j \in I} \tau^i_j $,   we get
\begin{equation}\label{eq27}
f_i(x)-f_i(x_0)+\epsilon_i
+\sum\limits_{ j\in I,j\not = i}[\bar{\tau}^i_j \hat M f_j(x)-\bar{\tau}^i_j\hat M f_j(x_0)
+\bar{\tau}^i_j \hat M \epsilon_j]
+ \sum\limits_{r \in L} \bar{\mu}^i_r g_r(x) \geq 0,
\end{equation}
for all $x\in \mathbb{R}^n$, where $\bar{\tau}^i_j= \frac{\tau^i_j}{\sum\limits_{j\in I} \tau^i_j }$ and  $\bar{\mu}^i_r= \frac{\mu^i_r}{\sum\limits_{j\in I} \tau^i_j  }$. In particular, for $x=x_0$, Inequality~(\ref{eq27}) gives
$\epsilon_i
+\sum\limits_{j \in I, j\not = i}\bar{\tau}^i_j \hat M \epsilon_j
+ \sum\limits_{r\in L} \bar{\mu}^i_r g_r(x_0) \geq 0.$
By setting $\bar{\epsilon_i}=\epsilon_i
+\sum\limits_{j \in I, j\not = i}\bar{\tau}^i_j \hat M \epsilon_j
$, we get Part~($ii$) as 
$\sum\limits_{r \in L} \bar{\mu}^i_r g_r(x_0) \geq -\bar{\epsilon_i}$.  Further, 
Inequality~(\ref{eq27}) reduces to, for all $x\in \mathbb{R}^n$,
\begin{eqnarray}\label{c2eq28}
f_i(x)
+\sum\limits_{j \in I, j\not = i}\bar{\tau}^i_j \hat M f_j(x)
+ \sum\limits_{r\in L} \bar{\mu}^i_r g_r(x)+\bar{\epsilon_i} \geq f_i(x_0)+\sum\limits_{j \in I, j\not = i}\bar{\tau}^i_j \hat M f_j(x_0).
\end{eqnarray}
As $x_0$ is feasible to CMOP,  $\sum\limits_{r\in L} \bar{\mu}^i_r g_r(x_0) \leq 0$. Thus, Inequality~(\ref{c2eq28}) becomes
\begin{eqnarray*}
f_i(x)
+\sum\limits_{j \in I, j\not = i}\bar{\tau}^i_j \hat M f_j(x)
+ \sum\limits_{r\in L} \bar{\mu}^i_r g_r(x)+\bar{\epsilon_i} \geq f_i(x_0)+\sum\limits_{j \in I, j\not = i}\bar{\tau}^i_j \hat M f_j(x_0)+\sum\limits_{r\in L} \bar{\mu}^i_r g_r(x_0),
\end{eqnarray*}
which implies that for each $i\in I$ and for all $x\in \mathbb{R}^n$,
\begin{eqnarray}\label{c2eq29}
L_i^{\hat{M}}(x,\bar{\tau}^i,\bar{\mu}^i)+\bar{\epsilon_i} \geq L_i^{\hat{M}}(x_0,\bar{\tau}^i,\bar{\mu}^i).
\end{eqnarray}
Further, from Equation~\eqref{eq410}, we observe that for all $i\in I$ and any $\mu \in \mathbb{R}^l_+$
\begin{eqnarray*}
L_i^{\hat{M}}(x_0,\bar{\tau}^i,\mu )\leq f(x_0)+\sum\limits_{j \in I, j\not = i}\bar{\tau}^i_j\hat{M}f_j(x_0),
\end{eqnarray*}
which can be written as $L_i^{\hat{M}}(x_0,\bar{\tau}^i,\mu)\leq f(x_0)+\sum\limits_{j \in I, j\not = i}\bar{\tau}^i_j\hat{M}f_j(x_0)+ \sum\limits_{r\in L} \bar{\mu}^i_r g_r(x)+\bar \epsilon_i .$
Thus, for all $x\in \mathbb{R}^n$ and $\mu \in \mathbb{R}^l_+$,
\begin{eqnarray}\label{c2eq30}
L_i^{\hat{M}}(x_0,\bar{\tau}^i,\mu ) \leq L_i^{\hat{M}}(x_0,\bar{\tau}^i,\bar{\mu}^i)+\bar{\epsilon_i}.
\end{eqnarray}
The Inequalities~(\ref{c2eq29}) and (\ref{c2eq30}) together prove Part~($i$). Now, for the sufficient part, let us assume that for a given $x_0\in \mathbb{R}^n$ and each $i \in I$ there exists $\bar{\tau}^i\in S^m$ and $\bar{\mu}^i\in \mathbb{R}^l_+$ such that Conditions ($i$) and ($ii$) hold. Our first step is to show that $x_0$ is feasible to CMOP. As we know from $(i)$, for all $\mu \in \mathbb{R}^l_+$
\begin{eqnarray*}
L_i^{\hat{M}}(x_0,\bar{\tau}^i,\mu )-\bar{\epsilon_i} \leq L_i^{\hat{M}}(x_0,\bar{\tau}^i,\bar{\mu}^i).
\end{eqnarray*}
Thus,  $f_i(x_0)
+\sum\limits_{j \in I, j\not = i}\bar{\tau}^i_j \hat M f_j(x_0)
+ \sum\limits_{r\in L} {\mu}_r g_r(x_0)-\bar{\epsilon_i }\leq f_i(x_0)+\sum\limits_{j 
\in I, j\not = i}\bar{\tau}^i_j \hat M f_j(x_0).$
This shows that for all $ \mu\in \mathbb{R}^l_+$,
\begin{equation}\label{eq049}
\sum\limits_{r\in L} \mu_r g_r(x_0) \leq \bar{\epsilon_i}.
\end{equation}
On the contrary, suppose $x_0$ is not feasible. Then,  there exists $r_0\in L$  such that $g_{r_0}(x_0)>0$. Then, choose $\mu=(0,\ldots,0,\mu_{r_0},0,\ldots,0)$, with $\mu_{r_0}>0$ and sufficiently large such that  $\mu_{r_0}g_{r_0}(x_0)> \bar{\epsilon_i}.$ Note that this contradicts
 Inequality~\eqref{eq049}. Hence, we conclude that $x_0$ is a feasible solution of CMOP.

  Now from right hand side of ($i$) we also have, for all $x\in \mathbb{R}^n$
\begin{eqnarray}
L_i^{\hat{M}}(x,\bar{\tau}^i,\bar{\mu}^i)+\bar{\epsilon_i} \geq L_i^{\hat{M}}(x_0,\bar{\tau}^i,\bar{\mu}^i).
\end{eqnarray}
which implies
\begin{multline*}
f_i(x)
+\sum\limits_{j \in I, j\not = i}\bar{\tau}^i_j \hat M f_j(x)
+ \sum\limits_{r\in L} \bar{\mu}^i_r g_r(x)+\epsilon_i
+\sum\limits_{j \in I, j\not = i}\bar{\tau}^i_j \hat M \epsilon_j \geq f_i(x_0) \\
+\sum\limits_{j \in I, j\not = i}\bar{\tau}^i_j \hat M f_j(x_0)
+\sum\limits_{r\in L} \bar{\mu}^i_r g_r(x_0).
\end{multline*}
Now, for any feasible $x$, $\sum\limits_{r\in L} \bar{\mu}^i_r g_r(x)\leq 0$. Thus, from the above inequality we have,
\begin{eqnarray}\label{eq55}
f_i(x)
+\sum\limits_{j \in I, j\not = i}\bar{\tau}^i_j \hat M f_j(x)
+\epsilon_i
+\sum\limits_{j \in I, j\not = i}\bar{\tau}^i_j \hat M \epsilon_j \geq f_i(x_0) & \nonumber \\
& \hspace*{-1.5in} +\sum\limits_{j \in I, j\not = i}\bar{\tau}^i_j \hat M f_j(x_0)
+\sum\limits_{r\in L} \bar{\mu}^i_r g_r(x_0).
\end{eqnarray}
Using Condition ($ii$), we have 
\begin{multline*}
f_i(x)
+\sum\limits_{j \in I, j\not = i}\bar{\tau}^i_j \hat M f_j(x)
+\epsilon_i
+\sum\limits_{j \in I, j\not = i}\bar{\tau}^i_j \hat M \epsilon_j \geq f_i(x_0) \\ +\sum\limits_{j \in I, j\not = i}\bar{\tau}^i_j \hat M f_j(x_0)-(\epsilon_i
+\sum\limits_{j \in I, j\not = i}\bar{\tau}^i_j \hat M \epsilon_j).
\end{multline*}
Since, it holds for each $i$, by summing over all the $i$'s we get,
\begin{eqnarray*}
\sum\limits_{i\in I} (1
+\hat M \sum\limits_{j \in I, j\not = i}\bar{\tau}^i_j) f_i(x)
+\sum\limits_{i\in I} (1+ \hat M\sum\limits_{j \in I, j\not = i}\bar{\tau}^i_j) (2\epsilon_j) \geq \sum\limits_{i\in I} (1+\hat M\sum\limits_{j \in I, j\not = i}\bar{\tau}^i_j ) f_i(x_0).
\end{eqnarray*}
Hence, $x_0$ is $\langle s, 2\epsilon\rangle$-minimizer of $P(s)$,  where  $s=(s_1, \ldots, s_m)$ with $s_i=1
+\hat M \sum\limits_{k \in I,k\not = i} \bar{\tau}^i_k$, for  $i \in I$. Now since $\bar{\tau}^i\in S^m$ for all $i,$ we have for all $i,j\in I$
$$\frac{s_i}{s_j}=\frac{1
+\hat M \sum\limits_{k\in I,k\not = i}\bar{\tau}^i_k}{1
+\hat M \sum\limits_{k\in I,k\not = j}\bar{\tau}^j_k}=\frac{1
+\hat M (1-\bar{\tau}^i_i)}{1
+\hat M (1-\bar{\tau}^j_j)}\leq 1+\hat{M}.$$
Since the above inequality is true for every $i$ and $j$, we have $\max\limits_{i,j}\{\frac{s_i}{s_j}\}\leq 1+\hat{M}$. Now consider $\tilde{M}\geq(1+\hat{M})(m-1)$ and using Theorem \ref{t2}, we conclude that 
$x_0\in \mathcal{G}_{\tilde{M},2\epsilon}(f,X).$ This completes the proof.\hfill $\Box$
\begin{remark}\rm
The saddle point type conditions are useful as a sufficient condition  if the number of objectives are only few in number. In fact, for sufficiency we can have a much simpler condition which we now state.
Let $x_0\in \mathbb{R}^n$ be a point that satisfies: \\ for each $i\in I$, there exists $\bar{\tau}^i\in S^m$ and $\bar{\mu}^i\in \mathbb{R}^l_+$ such that for all  $\mu \in \mathbb{R}^l_+ $ and $x\in \mathbb{R}^n$,
\begin{description}
\item $(a)$ $L_i^{\hat{M}}(x_0,\bar{\tau}^i,\mu)-\epsilon_i\leq L_i^{\hat{M}}(x_0,\bar{\tau}^i,\bar{\mu}^i)\leq L_i^{\hat{M}}(x,\bar{\tau}^i,\bar{\mu}^i)+\epsilon_i,$
\item $(b)$ $\sum\limits_{r\in L} \bar{\mu}^i_r g_r(x_0) \geq -{\epsilon_i}$.
\end{description}
Then, $x_0\in \mathcal{G}_{\tilde{M},2\epsilon}(f,X).$ 

  In order to prove the above statement, note that  $\bar{\epsilon_i}=\epsilon_i
+\sum\limits_{j=1,j\not = i}\bar{\tau}^i_j \hat M \epsilon_j
$. So, $\bar{\epsilon_i}\geq \epsilon_i$. Hence,  Conditions~($a$) and~($b$) above implies that  Conditions ($i$) and ($ii$) of Theorem \ref{t041}  are satisfied. Therefore,  we can simply apply the converse part of Theorem~\ref{t041} to get  $x_0\in \mathcal{G}_{\tilde{M},2\bar{\epsilon}}(f,X),$ where $\tilde{M}\geq (1+\hat{M})(m-1)$.  Note that  Condition~($a$) and ($b$) above are much simpler as  compared to checking Conditions~($i$) and ($ii$) as  $\bar{\epsilon_i}$ involves the multipliers $\bar{\tau}^i_j$. Hence, for the sufficiency part of  Theorem~\ref{t041} which requires the verification of Conditions~$(i)$ and ($ii$),  we will be using Conditions~($a$) and ($b$). 

Of course from the necessary part of Theorem \ref{t041}, we can also derive a multiplier rule involving $\epsilon$-subdifferentials, however this rule will be quite different. Observe that if $x_0\in \mathcal{G}_{\hat{M},\epsilon}(f,X)$, then Condition~($i$) of Theorem~\ref{t041} implies that for any $i\in I$ there exists $ \bar{\tau}^i\in S^m$ and $\bar{\mu}^i\in \mathbb{R}^l_+$ such that for all $x\in \mathbb{R}^n$,
$$ L_i^{\hat{M}}(x_0,\bar{\tau}^i,\bar{\mu}^i)\leq L_i^{\hat{M}}(x,\bar{\tau}^i,\bar{\mu}^i)+\bar{\epsilon_i},$$
which implies that $x_0\in \bar{\epsilon_i}- \underset{x\in \mathbb{R}^n}{\arg\min}L^{\hat{M}}_i(\; \cdot ,\bar{\tau}^i,\bar{\mu}^i),$ where $\bar \epsilon_i-\arg\min$ is the set of $\bar \epsilon_i$-minima of the function $L_i^{\hat{M}}(x,\bar{\tau}^i,\bar{\mu}^i)$. Thus, for each $i \in I$, $0\in \partial_{\bar{\epsilon_i}}L_i^{\hat{M}}(x_0,\bar{\tau}^i,\bar{\mu}^i).$
In fact a more compact necessary condition of the KKT type is given as follows, 
\begin{eqnarray}\label{eq33}
0\in \sum_{i=\in I}\partial_{\bar{\epsilon_i}}L_i^{\hat{M}}(x_0,\bar{\tau}^i,\bar{\mu}^i)\quad 
\text{with} \quad
\sum\limits_{r\in L} \bar{\mu}^i_r g_r(x_0) \geq -\bar{\epsilon_i}.
\end{eqnarray}
\end{remark}

\begin{thm}\label{t42}
For a given $\epsilon \in \mathbb{R}^m_+$ and $\hat{M}>0$, let us consider the problem CMOP. If $x_0\in \mathcal{G}_{\hat{M},{\epsilon}}(f,X)$, then there exist vectors  $\bar{\tau}^i\in S^m$ and $\bar{\mu}^i\in \mathbb{R}^l_+$, $i\in I$ such that 
\begin{description}
\item $(A)$ $0\in \sum\limits_{i\in I}\partial_{\bar{\epsilon_i}}L_i^{\hat{M}}(x_0,\bar{\tau}^i,\bar{\mu}^i),$
\item $(B)$ $\sum\limits_{r=1}^l \bar{\mu}^i_r g_r(x_0) \geq -\bar{\epsilon_i}$,
\end{description}
where $\bar{\epsilon_i}=\epsilon_i
+\sum\limits_{j\in I,j\not = i}\bar{\tau}^i_j \hat M \epsilon_j
$, $i\in I$. Conversely, if $x_0\in X$ be a point for which  there exist vectors $(\bar{\tau}^i,\bar{\mu}^i)\in S^m\times \mathbb{R}^l_+$, $i\in I$ such that ($A$) and ($B$) hold then $x_0\in \mathcal{G}_{\tilde{M},{2\epsilon}}(f,X)$, where $\tilde{M}=(1+\hat{M})(m-1)$.
\end{thm}
\textit{Proof:}
The necessary part has already been done in above remark. For sufficient part, let conditions $(A)$ and $(B)$ hold for $x_0\in X$. This means that there exists $ \bar{v}^i\in\partial_{\bar{\epsilon_i}}L_i^{\hat{M}}(x_0,\bar{\tau}^i,\bar{\mu}^i) $ for all $i\in I$ such that 
\begin{equation}\label{eq35}
0=\bar{v}^1+\bar{v}^2+\ldots+\bar{v}^m.
\end{equation}
Thus, from definition of $\epsilon$-subdifferential,  for each $i\in I$,
 $$L^{\hat{M}}_i(x,\bar{\tau},\bar{\mu}^i)-L^{\hat{M}}_i(x,\bar{\tau}^i,\bar{\mu}^i)\geq \langle\bar{v}^i,x-x_0\rangle-\bar{\epsilon}^i.$$
 Hence, 
 $$\sum\limits_{i\in I} L^{\hat{M}}_i(x,\bar{\tau},\bar{\mu}^i)-\sum\limits_{i \in I} L^{\hat{M}}_i(x,\bar{\tau}^i,\bar{\mu}^i)\geq \langle \sum\limits_{i \in I}\bar{v}^i,x-x_0\rangle-\sum\limits_{i\in I}\bar{\epsilon}^i.$$
Now using Equation~(\ref{eq35}), we get 
\begin{multline*}
\sum\limits_{i \in I} (f_i(x)
+\sum\limits_{j \in I, j\not = i}\bar{\tau}^i_j \hat M f_j(x)
+ \sum\limits_{r \in L} \bar{\mu}^i_r g_r(x))
- \sum\limits_{i \in I} (f_i(x_0) \\
+\sum\limits_{j \in I, j\not = i}\bar{\tau}^i_j \hat M f_j(x_0) 
+\sum\limits_{r\in L} \bar{\mu}^i_r g_r(x_0))\geq -\sum\limits_{i \in I} \bar{\epsilon}_i.
\end{multline*}
So, if $x$ is a feasible point then using Condition~($B$), the above inequality reduces to
\begin{eqnarray*}
\sum\limits_{i \in I} (f_i(x)
+\sum\limits_{j \in I, j\not = i}\bar{\tau}^i_j \hat M f_j(x))\geq 
\sum\limits_{i \in I} (f_i(x_0)
+\sum\limits_{j\in I, j\not = i}\bar{\tau}^i_j \hat M f_j(x_0))-\sum\limits_{i \in I} 2\bar{\epsilon}_i,
\end{eqnarray*} 
which can be rewritten as
\begin{eqnarray*}
\sum\limits_{i \in I} (1
+\sum\limits_{i \in I, i\not = j}\bar{\tau}^i_j \hat M) f_i(x)\geq 
\sum\limits_{i \in I} (1
+\sum\limits_{i \in I, i\not = j}\bar{\tau}^i_j \hat M )f_i(x_0)-\sum\limits_{i \in I} (1
+\sum\limits_{i \in I, i\not = j}\bar{\tau}^i_j \hat M) 2\epsilon_i.
\end{eqnarray*} 
Hence, $x_0$ is $\langle s, 2\epsilon\rangle$-minimizer of $P(s)$ where $s_i=1
+\hat M \sum\limits_{k\in I,k\not = i} \bar{\tau}^i_k$. Now using the same argument as in Theorem~\ref{t2}, we conclude that $x_0\in \mathcal{G}_{\tilde{M},{2\epsilon}}(f,X)$, where $\tilde{M}=(1+\hat{M})(m-1)$. This completes the proof.\hfill$\Box$

\section{Concluding remarks}

To analyze the behaviour of an optimization problem from the viewpoint of KKT conditions is deep-rooted in psyche of researchers in optimization theory. Though KKT conditions may not have been used very heavily in multiobjective optimization, but they can, however, act very well as a tool to develop stopping criteria. In this article, we characterize approximate versions of Pareto and proper Pareto solution using KKT type conditions. In fact, in the convex case, we achieve a complete characterization, for example, Theorem \ref{t22} demonstrates that a sequence of points which converge to weak Pareto minimizer has a subsequence where each point satisfies an approximate version of the KKT conditions. This result thus demonstrates the reason why approximate KKT type conditions can be used as stopping criteria.\\
The analysis of the approximate versions of the $\hat{M}$-Geoffrion proper solutions in terms of approximate KKT conditions is a starting point for building stopping criteria to identify such points. Our future research would involve more computational studies by using these optimality conditions as a stopping criterion.

\end{document}